\documentclass{amsart}
\usepackage{amsfonts}

\usepackage{graphicx}
\usepackage{amscd}
\usepackage{amsmath}

\newtheorem{theorem}{Theorem}
\theoremstyle{plain}

\newtheorem{corollary}[theorem]{Corollary}

\newtheorem{fact}[theorem]{Fact}

\theoremstyle{definition}
\newtheorem{remark}[theorem]{Remark}

\newcounter{thmenumerate}
\begin{document}

\title[Covering numbers in Schauder's Theorem about adjoints]{Estimates for covering numbers in Schauder's theorem about adjoints
of compact operators}

\author{Michael Cwikel and Eliahu Levy}

\address{Department of Mathematics, Technion - Israel Institute of Technology,
Haifa 32000, Israel}

\email{mcwikel@math.technion.ac.il, eliahu@techunix.technion.ac.il}

\thanks{The research of the first named author was supported by the Technion
V.P.R.\ Fund and by the Fund for Promotion of Research at the Technion. }



\maketitle
\begin{abstract}
Let $T:X\to Y$ be a bounded linear map between Banach spaces $X$ and $Y$.
Let $T^*:Y^* \to X^*$ be its adjoint. Let $\mathcal{B}_X$ and $\mathcal{B}_{Y^*}$
be the closed
unit balls of $X$ and $Y^*$ respectively.
We obtain apparently new estimates for the covering numbers of
the set $T^*\left(\mathcal{B}_{Y^*}\right)$.
These are expressed in terms of the covering numbers
of $T\left(\mathcal{B}_X\right)$, or, more generally, in terms of the covering numbers of a
``significant" subset of $T\left(\mathcal{B}_X\right)$.
The latter more general estimates are best possible. These estimates
follow from our new quantitative version of an abstract compactness
result which generalizes classical theorems of Arzel\`a-Ascoli and of
Schauder. Analogous estimates also hold for the covering numbers of
$T\left(\mathcal{B}_X\right)$, in terms of the covering numbers of
$T^*\left(\mathcal{B}_{Y^*}\right)$
or in terms of a
suitable ``significant" subset of $T^*\left(\mathcal{B}_{Y^*}\right)$.
\end{abstract}

\section{\label{sec:intro}Introduction}

The main motivation for the main result of this note is to give quantitative
versions of the celebrated Schauder theorem about adjoints of compact
operators. In fact our versions also apply to operators which are
not compact.

When we first obtained these results it seemed hard to imagine that
they are not already known. But so far we have not found any references
to similar results in the literature. We invite the reader to inform
us of any such references. In future versions of this note we hope
to include a more extensive survey of related previous results.

Of course there are remarkable recent results by Shiri Artstein--Avidan,
Vitali Milman, Stanislaw Szarek and Nicole Tomczak--Jaegerman \cite{SMS,SMST}
which give estimates for covering numbers which look much much stronger
than the ones we shall give here. But they deal with a slightly different
kind of question. Still more recently, Emanuel Milman \cite{emanuel}
has obtained results which in many cases also give substantially better
estimates than those given below in our Corollary \ref{cor:schauder}.
However his estimates contain factors and exponents which depend on
the dimension of the underlying space, and ours do not, at least not
in any explicit way. Moreover, we obtain a quantitative version of
a natural variant of Schauder's Theorem, (see Corollary \ref{cor:qschauder})
where our estimates are best possible.

Our point of departure is a theorem which is general enough to contain
the classical theorems of Arzelà-Ascoli and of Schauder as special
cases. But it, in turn, can be considered as a special case of considerably
more abstract and general results presented by Robert G.~Bartle in
\cite{bartle}, and which have their roots in earlier work of Kakutani,
R.\ S.\ Phillips and \v{S}mulian.

The said theorem appears in \cite{elevy} as a prelude to other results
dealing with finitely additive means and semigroups of operators.
It also appears in \cite{clat}, where it is used as a tool to study
complex interpolation of compact operators. Each of us obtained the
theorem independently, and before we became aware of the earlier and
more general results of \cite{bartle}.

We will state the theorem in the same formulation as is used in \cite{clat}.
Perhaps we first need to recall that a \textit{semimetric space} $(E,d)$
(also often referred to as a \textit{pseudometric space}) is a set
$E$ equipped with a function $d:E\times E\to[0,\infty)$ which is
a \textit{semimetric}, this meaning that it satisfies all the usual
conditions for a metric, except that the condition $d(x,y)=0$ does
not imply that $x=y$. The definition of a totally bounded semimetric
space (which will not be explicitly needed anyway in the sequel) is
exactly analogous to that of a totally bounded metric space.

\smallskip{}

\begin{theorem}
\label{aas} (Cf.~\cite{clat,elevy}) Let $A$ and $B$ be two sets
and let $h:A\times B\rightarrow\mathbb{C}$ be a function with the
properties that \begin{equation}
\sup_{a\in A}\left|h(a,b)\right|<\infty\text{ for each fixed }b\in B\text{, and}\label{qzi}\end{equation}
 \begin{equation}
\sup_{b\in B}\left|h(a,b)\right|<\infty\text{ for each fixed }a\in A.\label{qzii}\end{equation}

\noindent Define $d_{A}(a_{1},a_{2}):=\sup_{b\in B}\left|h(a_{1},b)-h(a_{2},b)\right|$
for each pair of elements $a_{1}$ and $a_{2}$ in $A$.

\noindent Define $d_{B}(b_{1},b_{2})=\sup_{a\in A}\left|h(a,b_{1})-h(a,b_{2})\right|$
for each pair of elements $b_{1}$ and $b_{2}$ in $B$.

Then \begin{equation}
(A,d_{A})\text{ and }(B,d_{B})\text{ are semimetric spaces}\label{qziii}\end{equation}

\noindent and \begin{equation}
(A,d_{A})\text{ is totally bounded if and only if }(B,d_{B})\text{ is totally bounded. }\label{qziv}\end{equation}

\end{theorem}
The claim (\ref{qziii}) in the previous theorem is of course a trivial
consequence of the definitions of $d_{A}$ and $d_{B}$ and we shall
use it in the sequel here without further comment.

There are proofs of Theorem \ref{aas} in both \cite{clat} and \cite{elevy}.
The advantage of the proof in \cite{elevy} is that it can be adapted
to give a quantitative version of the theorem, and therefore also
of Schauder's Theorem. That is what we do in this note.

Before we can state our new results we need to fix and discuss some
more (essentially standard) notation and terminology:

For each Banach space $X$ we denote the closed unit ball of $X$
by $\mathcal{B}_{X}$.

We use the standard notation of the ``\textit{ceiling function}'',
i.e., for each $t\in\mathbb{R}$ we let $\left\lceil t\right\rceil $
denote the smallest integer which dominates $t$.

Let $(E,d)$ be a semimetric space. We wish to define balls in $(E,d)$,
the diameters of subsets of $E$, and two kinds of covering numbers
for $E$. Our definitions will be mostly obvious variants of familiar
ones for metric spaces. Covering numbers in metric spaces are often
defined in terms of coverings by open balls of fixed given radius,
(whose centres may or may not be required to be in some particular
subset being considered). Here we will find it convenient to instead
use coverings by ``closed balls'' or arbitrary sets of fixed given
diameter. We have permitted ourselves these slight ``perturbations''
of the usual definitions in order to enable the convenient formulation
of examples showing that at least some of our results are best possible.
Of course our results can easily be translated into results corresponding
to covering numbers defined via coverings using open balls.

For each subset $G$ of $E$ we define the \textit{diameter} of $G$
to (of course!) be the quantity\[
\mathrm{diam}(G)=\sup_{x,y\in G}d(x,y)\,.\]

For each $x\in E$ and each $r>0$ we refer to the sets $\left\{ y\in E:d(x,y)<r\right\} $
and $\left\{ y\in E:d(x,y)\le r\right\} $ respectively as the \textit{open}
and \textit{closed balls of radius $r$ centred at }$x$.

For each semimetric space $(E,d)$ and each $\epsilon>0$, the \textit{intrinsic
covering number} $N_{E}(\epsilon)$ is the least positive integer
$n$ for which there exists a finite subset $F\subset E$ of cardinality
$n$ such that \begin{equation}
\min_{y\in F}d(x,y)\le\epsilon\mbox{ for each }x\in E\,,\label{eq:hb}\end{equation}
i.e., $N_{E}(\epsilon)$ is the smallest $n$ such that $E$ is contained
in some union of $n$ closed balls of radius $\epsilon$. If no such
$n\in\mathbb{N}$ exists, then $N_{E}(\epsilon)=\infty$.

We use the word ``intrinsic'' in this definition to reflect the fact
that if $E$ happens to be contained in a larger semimetric space
and $d$ is the restriction to $E\times E$ of the semimetric for
that larger space, then we are requiring the centres of all the closed
balls to be in $E$. This requirement will be relevant in the settings
of Corollaries \ref{cor:schauder} and \ref{cor:qschauder}.

For each semimetric space $(E,d)$ and each $\epsilon>0$ we define
the \textit{diameter covering number} $N_{E}^{\Delta}(\epsilon)$
to be the smallest positive integer $n$ for which there exist $n$
subsets $E_{1}$, $E_{2}$, ....., $E_{n}$ of $E$, each having diameter
not exceeding $2\epsilon$ and for which $E\subset\bigcup_{j=1}^{n}E_{j}$.
If no such positive integer $n$ exists, then $N_{E}^{\Delta}(\epsilon)=\infty$.

The triangle inequality for the semimetric $d$ immediately gives
us some simple connections between intrinsic covering numbers and
diameter covering numbers, namely

\begin{equation}
N_{E}(2\epsilon)\le N_{E}^{\Delta}(\epsilon)\le N_{E}(\epsilon)\,\mbox{for all }\epsilon>0\,.\label{eq:tkcn}\end{equation}

The naive hope that each arbitrary subset of $E$ having diameter
not exceeding $2\epsilon$ might be contained in a ball of radius
$\epsilon$ is immediately shattered by the following example: In
$\mathbb{R}^{2}$ with euclidean norm take $E$ to be the interior
of an equilateral triangle of side length $2\epsilon$. However we
can make the following simple observation, which will turn out to
be relevant later in the setting where we will show that one of our
results is best possible.

\begin{fact}
\label{fac:texz}Let $m$ be a positive integer, let $d(x,y)=\left\Vert x-y\right\Vert _{\ell_{m}^{\infty}}$
for each $x,y\in\mathbb{R}^{m}$ and let $E$ be a subset of $\mathbb{R}^{m}$.
Let $G$ be a subset of $E$ with $\mathrm{diam}(G)\le2\epsilon$.
Then there exists a point $x\in\mathbb{R}^{m}$ such that $G\subset\left\{ y\in E:d(x,y)\le\epsilon\right\} $.
\end{fact}
\textit{Proof.} For each $k=1,2,....,m$ and each $x\in\mathbb{R}^{m},$
let $\pi_{k}(x)$ denote the $k^{th}$ coordinate of $x$. We let
$\alpha_{k}=\inf\left\{ \pi_{k}(x):x\in G\right\} $ and $\beta_{k}=\sup\left\{ \pi_{k}(x):x\in G\right\} $.
The finiteness of $\mathrm{diam}(G)$ ensures that both $\alpha_{k}$
and $\beta_{k}$ are finite. Then the bound on $\mathrm{diam}(G)$
ensures that $0\le\beta_{k}-\alpha_{k}\le2\epsilon$. The point $x=\left(\frac{\alpha_{1}+\beta_{1}}{2},\frac{\alpha_{2}+\beta_{2}}{2},.....,\frac{\alpha_{m}+\beta_{m}}{2}\right)$
clearly has the required property. $\qed$

\begin{remark}
Note that the point $x$ in Fact \ref{fac:texz} does not have to
belong to $G$ or even to $E$. If $G$ happens to be convex and closed
then we can have $x\in G$ for $m\le2$. But this need not happen
for $m\ge3$ as one can see by considering the set $G\subset\mathbb{R}^{3}$
which is the convex hull of $\left(-1,1,1\right)$, $\left(1,-1,1\right)$
and $(1,1,-1)$.
\end{remark}
\bigskip{}

We are now ready to state our main result, and then its corollaries.
Their proofs will be given in Section \ref{sec:proofs}.

\begin{theorem}
\label{thm:elq}Let $A$ and $B$ be two sets and let $h:A\times B\rightarrow\mathbb{C}$
be a function with the properties stated in Theorem \ref{aas}. Let
$d_{A}$ and $d_{B}$ be the semimetrics defined on $A$ and $B$
respectively, as in Theorem \ref{aas}.

Suppose that the intrinsic covering number $N_{A}(\epsilon)$ is finite
for some positive number $\epsilon$. Then

(i) The quantity $C:=\sup_{a\in A,b\in B}\left|h(a,b)\right|$ is
also finite.

(ii) The diameter covering number $N_{B}^{\Delta}(\rho)$ is finite
for each $\rho>\epsilon$.

(iii) Furthermore,

\begin{equation}
N_{B}^{\Delta}(\epsilon+\delta)\le\left(\left\lceil \frac{\sqrt{2}C}{\delta}\right\rceil \right)^{2N_{A}(\epsilon)}\,\mbox{for each }\delta>0\,,\label{eq:cgkt}\end{equation}
 and, if $h$ is real valued, the following stronger estimate also
holds. \begin{equation}
N_{B}^{\Delta}(\epsilon+\delta)\le\left(\left\lceil \frac{C}{\delta}\right\rceil \right)^{N_{A}(\epsilon)}\,\mbox{for each }\delta>0\,.\label{eq:gkt}\end{equation}

(iv) By symmetry, the roles of $A$ and $B$ can be interchanged and
so exactly analogous estimates hold for $N_{A}^{\Delta}(\epsilon+\delta)$
in terms of $N_{B}(\epsilon)$.
\end{theorem}

\begin{remark}
Note that in this theorem we do not make any ``compactness'' or ``total
boundedness'' assumptions about $(A,d_{A})$ or $(B,d_{B})$.
\end{remark}

\begin{remark}
\label{rem:absconv}In general, in Theorem \ref{thm:elq}, we cannot
expect to have any estimate for the supremum $C$ just in terms of
$N_{A}(\epsilon)$. However we do have an estimate for $C$ in the
following particular case: Assume that $A$ is an absolutely convex
subset of a linear space $V$ and that the semimetric $d_{A}$ defined
as above also satisfies \begin{equation}
d_{A}(a,a')=p(a-a')\label{eq:ffda}\end{equation}
for some seminorm $p$ on $V$ and all $a,a'\in A$. Assume furthermore
that the function $h:A\times B\to\mathbb{C}$ satisfies \begin{equation}
h(0,b)=0\,\,\mbox{for all }b\in B.\label{eq:hbz}\end{equation}
In this case, from straightforward ``geometric'' considerations,
one can expect the supremum $C$ to satisfy\begin{equation}
C\le\epsilon N_{A}(\epsilon)\,,\label{eq:zefc}\end{equation}
and we will prove in an appendix (Section \ref{sec:appendix}) that
(\ref{eq:zefc}) indeed holds. Consequently, in this case, Theorem
\ref{thm:elq} gives us that \begin{equation}
N_{B}^{\Delta}(\epsilon+\delta)\le\left(\left\lceil \frac{\sqrt{2}\epsilon N_{A}(\epsilon)}{\delta}\right\rceil \right)^{2N_{A}(\epsilon)}\,,\label{eq:cbi2}\end{equation}
and, if we also have that $h$ is real valued, we can sharpen this
to \begin{equation}
N_{B}^{\Delta}(\epsilon+\delta)\le\left(\left\lceil \frac{\epsilon N_{A}(\epsilon)}{\delta}\right\rceil \right)^{N_{A}(\epsilon)}.\label{eq:fvs}\end{equation}
In our first application (Corollary \ref{cor:schauder}) of Theorem
\ref{thm:elq} to the context of Schauder's Theorem, the set $A$,
the function $h$ and the semimetric $d_{A}$ satisfy the conditions
of this particular case. But of course there we have a simpler and,
in general, better way to estimate $C$. We apparently cannot invoke
an analogue of (\ref{eq:zefc}) in Corollary \ref{cor:qschauder}.
\end{remark}

\begin{remark}
One of the examples to be presented in Section \ref{sec:example}
will show that the estimate (\ref{eq:gkt}) cannot be sharpened. On
the other hand, it will be clear from the proof below that the estimates
(\ref{eq:cgkt}) and (\ref{eq:cbi2}) can be somewhat improved, either
by making some small efforts in planar geometry or by consulting some
appropriate references, in order to sharpen the claim (ii) in Fact
\ref{fac:disks}. Another example in Section \ref{sec:example} will
show that the requirement that $\rho>\epsilon$ in part (ii) of Theorem
\ref{thm:elq} cannot be replaced by any weaker requirement of the
form $\rho>\phi(\epsilon)$ for some function $\phi$ of $\epsilon$
which satisfies $\phi(\epsilon)<\epsilon$.
\end{remark}

Here is the obvious simplest way that we can apply Theorem \ref{thm:elq}.

\begin{corollary}
\label{cor:schauder}Let $X$ and $Y$ be Banach spaces and let $\mathcal{B}_{X}$
and $\mathcal{B}_{Y^{*}}$ be the closed unit balls of $X$ and $Y^{*}$
respectively.

Let $T:X\to Y$ be a bounded linear operator with adjoint $T^{*}:Y^{*}\to X^{*}$.

For each $\epsilon>0$ let $N_{T}(\epsilon)$ denote the least number
of closed balls in $Y$ of radius $\epsilon$ with centres in $T\left(\mathcal{B}_{X}\right)$
which are required to cover the set $T\left(\mathcal{B}_{X}\right)$,
and let $N_{T}^{\Delta}(\epsilon)$ denote the least number of subsets
of $Y$ each with $Y$-norm diameter not exceeding $2\epsilon$ which
are required to cover $T\left(\mathcal{B}_{X}\right)$.

Analogously, let $N_{T^{*}}(\epsilon)$ denote the least number of
closed balls in $X^{*}$ of radius $\epsilon$ with centres in $T^{*}\left(\mathcal{B}_{Y^{*}}\right)$
which are required to cover the set $T^{*}\left(\mathcal{B}_{Y^{*}}\right)$
and let $N_{T^{*}}^{\Delta}(\epsilon)$ denote the least number of
subsets of $X^{*}$ each with $X^{*}$-norm diameter not exceeding
$2\epsilon$ which are required to cover $T^{*}\left(\mathcal{B}_{Y^{*}}\right)$.

Suppose that $N_{T}(\epsilon)$ is finite for some particular $\epsilon>0$.
Then $N_{T^{*}}^{\Delta}(\rho)$ is finite for all $\rho>\epsilon$
and the estimate \begin{equation}
N_{T^{*}}^{\Delta}(\epsilon+\delta)\le\left(\left\lceil \frac{\sqrt{2}\left\Vert T\right\Vert _{X\to Y}}{\delta}\right\rceil \right)^{2N_{T}(\epsilon)}\label{eq:onz}\end{equation}
 holds for all $\delta>0$. If $X$ and $Y$ are real Banach spaces,
then this estimate can be sharpened to \begin{equation}
N_{T^{*}}^{\Delta}(\epsilon+\delta)\le\left(\left\lceil \frac{\left\Vert T\right\Vert _{X\to Y}}{\delta}\right\rceil \right)^{N_{T}(\epsilon)}\,.\label{eq:twz}\end{equation}

Furthermore, if $N_{T^{*}}(\epsilon)$ is finite for some $\epsilon>0$,
then $N_{T}^{\Delta}(\rho)$ is finite for all $\rho>\epsilon$ and
the quantity $N_{T}^{\Delta}(\epsilon+\delta)$ can be estimated in
terms of $N_{T^{*}}(\epsilon)$ via formulae exactly analogous to
(\ref{eq:onz}) and (\ref{eq:twz}), where $T$ and $T^{*}$ are interchanged.
\end{corollary}

Apparently other results will give much better estimates than (\ref{eq:onz})
and (\ref{eq:twz}). But here is a slightly more subtle variant of
Corollary \ref{cor:schauder} for which, in some cases, our estimates
are best possible. With the perspective of Theorem \ref{thm:elq}
we can see that it may be just as appropriate and just as easy to
work with the covering numbers of certain ``significant'' subsets
of $T\left(\mathcal{B}_{X}\right)$ and of $T^{*}\left(\mathcal{B}_{Y^{*}}\right)$,
instead of working with the covering numbers of these sets themselves.
We will obtain new versions of the estimates (\ref{eq:onz}) and (\ref{eq:twz})
for $N_{T^{*}}^{\Delta}(\epsilon+\delta)$, which are stronger in
the sense that the number $N_{T}(\epsilon)$ is replaced by a smaller,
in some cases very much smaller number, which is the covering number
of a suitable subset $K$ of $T\left(\mathcal{B}_{X}\right)$. Similarly
the estimates for $N_{T}^{\Delta}(\epsilon+\delta)$, which were stated
implicitly in Corollary \ref{cor:schauder}, can be replaced by stronger
results where $N_{T^{*}}(\epsilon)$ is replaced by the covering number
of a suitable subset $K^{*}$ of $T^{*}\left(\mathcal{B}_{Y^{*}}\right)$.

\begin{corollary}
\label{cor:qschauder} Let $X$, $Y$, $T$, $N_{T}^{\Delta}(\epsilon)$,
and $N_{T^{*}}^{\Delta}(\epsilon)$ all be as specified in the statement
of Corollary \ref{cor:schauder}. Let $K$ be a ``norming'' subset
of $T(\mathcal{B}_{X})$, i.e., a subset with the property that \begin{equation}
\sup\left\{ \left|\left\langle u,y\right\rangle \right|:u\in K\right\} =\sup\left\{ \left|\left\langle u,y\right\rangle \right|:u\in T(\mathcal{B}_{X})\right\} \mbox{ for each }y\in Y^{*}\,.\label{eq:npo}\end{equation}
 Analogously, let $K^{*}$ be a subset of $T^{*}\left(\mathcal{B}_{Y^{*}}\right)$
with the property that \begin{equation}
\sup\left\{ \left|\left\langle x,v\right\rangle \right|:v\in K^{*}\right\} =\sup\left\{ \left|\left\langle x,v\right\rangle \right|:v\in T^{*}(\mathcal{B}_{Y^{*}})\right\} \mbox{ for each }x\in X\,.\label{eq:npt}\end{equation}
For each $\epsilon>0$ let $N[K,\epsilon]$ be the least number of
closed balls in $Y$ of radius $\epsilon$ with centres in $K$ which
are required to cover the set $K$. Analogously, let $N[K^{*},\epsilon]$
denote the least number of closed balls in $X^{*}$ of radius $\epsilon$
with centres in $K^{*}$ which are required to cover the set $K^{*}$.

Suppose that $N[K,\epsilon]$ is finite for some particular $\epsilon>0$.
Then $N_{T^{*}}^{\Delta}(\rho)$ is finite for all $\rho>\epsilon$
and the estimate\begin{equation}
N_{T^{*}}^{\Delta}(\epsilon+\delta)\le\left(\left\lceil \frac{\sqrt{2}\left\Vert T\right\Vert _{X\to Y}}{\delta}\right\rceil \right)^{2N[K,\epsilon]}\,.\label{eq:qo}\end{equation}
holds for all $\delta>0$. If $X$ and $Y$ are real Banach spaces
then this estimate can be sharpened to\begin{equation}
N_{T^{*}}^{\Delta}(\epsilon+\delta)\le\left(\left\lceil \frac{\left\Vert T\right\Vert _{X\to Y}}{\delta}\right\rceil \right)^{N[K,\epsilon]}\,.\label{eq:qt}\end{equation}
 Analogously, if $N[K^{*},\epsilon]$ is finite for some particular
$\epsilon>0$, then $N_{T}^{\Delta}(\rho)$ is finite for all $\rho>\epsilon$
and the estimate\begin{equation}
N_{T}^{\Delta}(\epsilon+\delta)\le\left(\left\lceil \frac{\sqrt{2}\left\Vert T\right\Vert _{X\to Y}}{\delta}\right\rceil \right)^{2N[K^{*},\epsilon]}\label{eq:aqo}\end{equation}
 holds for all $\delta>0$, and in the case where $X$ and $Y$ are
real Banach spaces, it can be sharpened to \begin{equation}
N_{T}^{\Delta}(\epsilon+\delta)\le\left(\left\lceil \frac{\left\Vert T\right\Vert _{X\to Y}}{\delta}\right\rceil \right)^{N[K^{*},\epsilon]}\,.\label{eq:aqt}\end{equation}

\end{corollary}

\begin{remark}
Obviously Corollary \ref{cor:schauder} is nothing more than a special
case of Corollary \ref{cor:qschauder} since of course the sets $K=T\left(\mathcal{B}_{X}\right)$
and $K^{*}=T^{*}\left(\mathcal{B}_{Y^{*}}\right)$ satisfy (\ref{eq:npo})
and (\ref{eq:npt}). But it seems better and clearer to have begun
this discussion by stating that special case separately.
\end{remark}

What reasonable kinds of sets might play the roles of $K$ and of
$K^{*}$ in Corollary \ref{cor:qschauder}? Obviously we can take
$K=\left\{ Tx:\left\Vert x\right\Vert _{X}=1\right\} $ and $K^{*}=\left\{ T^{*}y:\left\Vert y\right\Vert _{Y^{*}}=1\right\} $.
But for such choices we cannot expect $N[K,\epsilon]$ and $N[K^{*},\epsilon]$
to be substantially smaller than, respectively, $N_{T}(\epsilon)$
and $N_{T^{*}}(\epsilon)$. In finite dimensional spaces, if $\mathcal{B}_{X}$
is the convex hull of some finite set $F$, then $K$ can be chosen
to be $T(F)$. (We are aware of at least two papers, namely \cite{ckp}
and \cite{kyrezi}, which could be applied to give connections between
the covering numbers, or entropy numbers, of $T(F)$ and of $T\left(\mathcal{B}_{X}\right)$
in such cases.) If $T$ has some special additional properties it
might be possible to choose an even smaller set than $T(F)$ in the
role of $K$.

Here is what is probably the most natural example of a choice of $K$
which satisfies (\ref{eq:npo}) and for which $N[K,\epsilon]$ is
very significantly smaller than $N_{T}(\epsilon)$. Let $X$ and $Y$
both be $\mathbb{R}^{n}$ equipped with the $\ell^{1}$ norm and let
$T$ be the identity operator on $\mathbb{R}^{n}$. Let $K$ be the
subset of $\mathcal{B}_{X}$ which consists of the $n$ points $e_{j}$
for $j=1,2,...,n$, where $e_{1}=(1,0,0,....,0)$, $e_{2}=(0,1,0,0,....,0)$,
....., $e_{n}=(0,0,....,0,1)$. Of course $N_{T}(\epsilon)$ is arbitrarily
large for small values of $\epsilon$. But $N[K,\epsilon]=n$ for
all $\epsilon$ in the range $0<\epsilon<1$. Of course in (\ref{eq:npo})
we take $\left\langle \cdot,\cdot\right\rangle $ to be the usual
inner product on $\mathbb{R}^{n}$, and so $X^{*}$ and $Y^{*}$ are
both $\mathbb{R}^{n}$ equipped with the $\ell^{\infty}$ norm. Clearly
(\ref{eq:npo}) holds here since, for each $y\in\mathbb{R}^{n}$,
both sides of (\ref{eq:npo}) equal $\left\Vert y\right\Vert _{\ell_{n}^{\infty}}$.

Let us now use this example to show that in general the estimate (\ref{eq:qt})
cannot be improved. In our context here (cf.~Fact \ref{fac:texz}),
for each $\rho>0$, we see that $N_{T^{*}}^{\Delta}(\rho)$ is the
minimal number of closed cubes in $\mathbb{R}^{n}$ of side length
$2\rho$ required to cover all of the closed cube $Q=\left[-1,1\right]^{n}$.
The interval $[-1,1]$ is contained in the union of $\left\lceil \frac{2}{2\rho}\right\rceil $
non overlapping closed intervals of length $2\rho$, and so it is
clear that $N_{T^{*}}^{\Delta}(\rho)\le\left(\left\lceil \frac{1}{\rho}\right\rceil \right)^{n}$.
By trivial considerations of volume, we must also have $(2\rho)^{n}N_{T^{*}}^{\Delta}(\rho)\ge2^{n}$.
It follows that, for each positive integer $m$, we have $N_{T^{*}}^{\mathbf{\Delta}}(\rho)=m^{n}$
for all numbers $\rho$ which satisfy $m^{n}-1<\rho^{-n}\le m^{n}$.
In particular this will hold whenever we choose $\rho=\left(m^{n}-\theta\right)^{-1/n}$
for some number $\theta\in(0,1)$. For our purposes we will choose
a particular value of $\theta\in(0,1)$ which is sufficiently small
to ensure that \begin{equation}
m-1<\left(m^{n}-\theta\right)^{1/n}<m\,.\label{eq:cte}\end{equation}
For each number $\epsilon$ which is in the range $0<\epsilon<\min\left\{ 1,\left(m^{n}-\theta\right)^{-1/n}\right\} $
we have $N[K,\epsilon]=n$ and we also have that the number $\delta=\delta(\epsilon)$
which satisfies $\epsilon+\delta=\left(m^{n}-\theta\right)^{-1/n}$
is positive. For each such number $\epsilon$ and for each corresponding
$\delta=\delta(\epsilon)$, the left side of (\ref{eq:qt}) equals
$m^{n}$ and the right side equals $\left(\left\lceil \frac{\left\Vert T\right\Vert _{X\to Y}}{\delta}\right\rceil \right)^{n}=\left(\left\lceil \frac{1}{\left(m^{n}-\theta\right)^{-1/n}-\epsilon}\right\rceil \right)^{n}=\left(\left\lceil \frac{(m^{n}-\theta)^{1/n}}{1-\epsilon(m^{n}-\theta)^{1/n}}\right\rceil \right)^{n}$.
Since we have chosen $\theta$ to satisfy (\ref{eq:cte}), we can
now choose $\epsilon$ sufficiently small, so that the number $\frac{(m^{n}-\theta)^{1/n}}{1-2\epsilon(m^{n}-\theta)^{1/n}}$
also lies in the open interval $(m-1,m)$. This makes the right side
of (\ref{eq:qt}) also equal to $m^{n}$. So indeed the estimate (\ref{eq:qt})
is best possible for certain values of the numbers $\epsilon$ and
$\delta$, and in fact for infinitely many such values, which can
be taken arbitrarily small.

\begin{remark}
If we change this example so that $X=Y$ is $\mathbb{R}^{n}$ equipped
with the $\ell^{\infty}$ norm, then we can of course take $K$ to
consist of the $2^{n}$ vertices of the cube of side length 2 centred
at the origin, i.e., the extreme points of $\mathcal{B}_{X}$. Or
$K$ may consist merely of half of these points. Perhaps we will look
at this example in more detail in a subsequent version of this paper.
\end{remark}

\section{\label{sec:proofs}Proofs }

\subsection{The proof of Corollary \ref{cor:qschauder} (and therefore also of
Corollary \ref{cor:schauder}) . \protect \\
}

It is rather obvious what needs to be done, even more so in the case
where $K=T\left(\mathcal{B}_{X}\right)$ and $K^{*}=T^{*}\left(\mathcal{B}_{Y^{*}}\right)$.
But let us write out the proof explicitly, at least in this preliminary
version of our paper:

Let $A=\left\{ x\in\mathcal{B}_{X}:Tx\in K\right\} $. We clearly
have $T(A)=K$. We also let $B=\mathcal{B}_{Y^{*}}$. Define $h:A\times B\to\mathbb{C}$
by $h(a,b)=\left\langle Ta,b\right\rangle $. Then we have \begin{equation}
\sup_{a\in A,b\in B}\left|h(a,b)\right|\le\left\Vert T\right\Vert _{X\to Y}=\left\Vert T^{*}\right\Vert _{Y^{*}\to X^{*}}\,,\label{eq:efc}\end{equation}
and also, by the Hahn--Banach Theorem, \begin{equation}
d_{A}(a_{1},a_{2})=\left\Vert Ta_{1}-Ta_{2}\right\Vert _{Y}\,\mbox{for all }a_{1},a_{2}\in A\,.\label{eq:qcmt}\end{equation}

Let $\epsilon$ be a positive number for which $N[K,\epsilon]$ is
finite and equals $n$. Then there exists a set of $n$ points $\left\{ y_{1},y_{2},...,y_{n}\right\} $
in $K$ such that $\min\left\{ \left\Vert y-y_{j}\right\Vert _{Y}:1\le j\le n\right\} \le\epsilon$
for each $y\in K$. Since $K=T\left(A\right)$, we have a set of $n$
points $\left\{ x_{1},x_{2},....,x_{n}\right\} $ in $A$ such that
$Tx_{j}=y_{j}$ for each $j=1,2,...,n$.

For each $a\in A$ we have $Ta\in K$ and therefore \[
\min\left\{ d_{A}\left(a,x_{j}\right):1\le j\le n\right\} =\min\left\{ \left\Vert Ta-Tx_{j}\right\Vert _{Y}:1\le j\le n\right\} \le\epsilon\,.\]
This shows that\begin{equation}
N_{A}(\epsilon)\le N[K,\epsilon]\,\mbox{for all }\epsilon>0\,.\label{eq:nat}\end{equation}

In order to obtain a formula for $d_{B}$, we use the fact that $h(a,b)$
is also equal to $\left\langle a,T^{*}b\right\rangle $. For each
pair of elements $b_{1},b_{2}\in B$ we again use the fact that $T(A)=K$
and we apply (\ref{eq:npo}) to obtain that\begin{eqnarray}
d_{B}(b_{1},b_{2}) & = & \sup_{a\in A}\left|\left\langle Ta,b_{1}-b_{2}\right\rangle \right|=\sup_{u\in K}\left|\left\langle u,b_{1}-b_{2}\right\rangle \right|=\sup_{u\in T\left(\mathcal{B}_{X}\right)}\left|\left\langle u,b_{1}-b_{2}\right\rangle \right|\nonumber \\
 & = & \sup_{x\in\mathcal{B}_{X}}\left|\left\langle Tx,b_{1}-b_{2}\right\rangle \right|=\sup_{x\in\mathcal{B}_{X}}\left|\left\langle x,T^{*}(b_{1}-b_{2})\right\rangle \right|=\left\Vert T^{*}b_{1}-T^{*}b_{2}\right\Vert _{X^{*}}\,.\label{eq:fvy}\end{eqnarray}

Next we will show that \begin{equation}
N_{T^{*}}^{\Delta}(\rho)\le N_{B}^{\Delta}(\rho)\,\mbox{for all}\,\rho>0\,.\label{eq:ret}\end{equation}
Let $\rho$ be an arbitrary positive number. If $m=N_{B}^{\Delta}(\rho)$
then there exists a collection $W_{1}$, $W_{2}$,...., $W_{m}$ of
$m$ subsets of $B=\mathcal{B}_{Y^{*}}$ such that $\sup_{b,b'\in W_{j}}d_{B}(b,b')\le2\rho$
for each $j=1,2,....,m$ and $B\subset\bigcup_{j=1}^{m}W_{j}$. Let
$V_{j}=T^{*}\left(W_{j}\right)$ for each $j$. Each set $V_{j}$
is contained in $T^{*}\left(\mathcal{B}_{Y^{*}}\right)$ and, in view
of (\ref{eq:fvy}), \begin{equation}
\sup_{v,v'\in V_{j}}\left\Vert v-v'\right\Vert _{X^{*}}=\sup_{b,b'\in W_{j}}\left\Vert T^{*}b-T^{*}b'\right\Vert _{X^{*}}=\sup_{b,b'\in W_{j}}d_{B}(b,b')\le2\rho\,.\label{eq:mbcd}\end{equation}
We also have $T^{*}\left(\mathcal{B}_{Y^{*}}\right)\subset\bigcup_{j=1}^{m}V_{j}$
which, together with (\ref{eq:mbcd}), gives us (\ref{eq:ret}).

Now to obtain (\ref{eq:qo}) and (\ref{eq:qt}) we simply have to
substitute (\ref{eq:ret}), (\ref{eq:efc}) and (\ref{eq:nat}) in
(\ref{eq:cgkt}) and (\ref{eq:gkt}) respectively.

The rest of the proof, i.e., of the estimates (\ref{eq:aqo}) and
of (\ref{eq:aqt}), is a variant of the arguments to obtain (\ref{eq:qo})
and (\ref{eq:qt}) that we have just presented. For this last part
of the proof we will choose $A=\left\{ y\in\mathcal{B}_{Y^{*}}:T^{*}y\in K^{*}\right\} $
and $B=\mathcal{B}_{X}$ and $h(a,b)=\left\langle Tb,a\right\rangle =\left\langle b,T^{*}a\right\rangle $,
and then use reasoning which is almost exactly analogous to that used
above. beginning with the observation that $T^{*}(A)=K^{*}$. We will
only need to make two small and obvious changes:

$\bullet$ The formula for $d_{A}$, i.e., $d_{A}(a_{1},a_{2})=\left\Vert T^{*}a_{1}-T^{*}a_{2}\right\Vert _{X^{*}}$,
whose analogue in (\ref{eq:qcmt}) followed from the Hahn--Banach
Theorem, here is simply a consequence of the definition of the norm
of a linear functional.

$\bullet$ The last step in (\ref{eq:fvy}) above simply used the
definition of the norm of a linear functional. The analogous step
here, which shows that $d_{B}\left(b_{1},b_{2}\right)=\left\Vert Tb_{1}-Tb_{2}\right\Vert _{Y}$,
will use the Hahn--Banach theorem.

We leave the details to the reader. $\qed$ \bigskip{}

\subsection{How to cover big intervals by small intervals, and big disks by small
disks.\protect \\
}

This very simple geometrical observation will be one of the components
in our proof of Theorem \ref{thm:elq}.

\begin{fact}
\label{fac:disks}Let $\delta$ and $C$ be numbers satisfying $0<\delta\le C$.

(i) The interval $[-C,C]$ can be contained in the union of $\left\lceil \frac{C}{\delta}\right\rceil $
closed intervals $I_{q}$, $q=1,2,....,\left\lceil \frac{C}{\delta}\right\rceil $,
each of length $2\delta$.

(ii) The closed disk $\left\{ z\in\mathbb{C}:\left|z\right|\le C\right\} $
can be contained in the union of $\left\lceil \frac{\sqrt{2}C}{\delta}\right\rceil ^{2}$
closed disks $D_{q}$, $q=1,2,....,\left\lceil \frac{\sqrt{2}C}{\delta}\right\rceil ^{2}$
each of radius $\delta$.
\end{fact}
\textit{Proof.} The claim (i) is trivial and needs no further comment.

The number of disks mentioned in the claim (ii) can certainly be reduced.
However, in this preliminary version of the paper, we will not seek
the optimal value of this quantity, but simply content ourselves with
some quite crude, simple minded estimates. First we remark that each
closed disc of radius $\delta$ contains a closed square of side length
$\sqrt{2}\delta$. If the positive integer $k$ satisfies $k\ge2C/\sqrt{2}\delta$
then we will be able to cover a closed square of side length $2C$
with $k^{2}$ closed squares of side length $\sqrt{2}\delta$. So
we will certainly be able to cover a closed disc of radius $C$ with
$k^{2}$ closed discs of radius $\delta$, where $k^{2}=\left\lceil \frac{\sqrt{2}C}{\delta}\right\rceil ^{2}$.
$\qed$

\subsection{The proof of Theorem \ref{thm:elq}. \protect \\
}

Throughout this proof the number $\epsilon>0$, for which we suppose
that $N_{A}(\epsilon)$ is finite, will remain fixed, and $n$ will
always equal $N_{A}(\epsilon)$. We also fix a finite $n$ element
subset $F=\left\{ a_{1},a_{2},.....,a_{n}\right\} $ of $A$ which
has the property \begin{equation}
\min\left\{ d_{A}(a,a_{j}):j\in\{1,2,...,n\}\right\} \le\epsilon\,.\label{eq:nfl}\end{equation}
At least one such subset necessarily exists.

$C$ will always denote the supremum $C:=\sup_{a\in A,b\in B}\left|h(a,b)\right|$.

Our first task is to show that the supremum $C$ is finite: Let $a$
and $a'$ be two arbitrary elements of $A$. Let $C_{F}=\max_{j,k\in\{1,2,....,n\}}d_{A}(a_{j},a_{k})$.
Then, for each $j$ and $k$ in $\{1,2,...,n\}$ we have \[
d_{A}(a,a')\le d_{A}(a,a_{j})+d_{A}(a_{j},a_{k})+d(a_{k},a')\,.\]
By making appropriate choices of $j$ and $k$ we obtain that \begin{equation}
d_{A}(a,a')\le2\epsilon+C_{F}\,.\label{eq:diama}\end{equation}

Now, for arbitrary $a\in A$ and $b\in B$, we have \begin{eqnarray*}
\left|h(a,b)\right| & \le & \left|h(a,b)-h(a_{1},b)\right|+|h(a_{1}b)|\\
 & \le & d_{A}(a,a_{1})+\sup_{b'\in B}\left|h(a_{1},b')\right|.\end{eqnarray*}
In view of (\ref{eq:diama}) and (\ref{qzii}) this last expression
is finite and does not depend on $a$ or $b$. So our assumption that
$N_{A}(\epsilon)$ is finite indeed ensures the finiteness also of
$C$.

Let $\delta$ be an arbitrary positive number. It remains to show
that $N_{B}^{\Delta}(\epsilon+\delta)$ is finite and satisfies the
required estimates.

It will be convenient to simultaneously consider both the general
case and the special case where $h$ is real valued. To facilitate
this we will let $\mathbb{K}$ denote the real field $\mathbb{R}$
if $h$ is real valued, and otherwise we will have $\mathbb{K}=\mathbb{C}$.

Let us now consider the set $G=\left\{ g_{b}:b\in B\right\} $ of
functions $g_{b}:F\to\mathbb{\mathbb{K}}$ defined by the formula
$g_{b}(f)=h(f,b)$. $G$ is of course contained in the set $Q\subset\ell^{\infty}(F)$
consisting of those functions $g:F\to\mathbb{K}$ for which $\left|g(f)\right|\le C$
for each $f$.

If $\mathbb{K}=\mathbb{C}$ then we let $M=\left\lceil \frac{\sqrt{2}C}{\delta}\right\rceil ^{2}$
and let $D_{1}$, $D_{2}$,....., $D_{M}$ be the $M$ closed discs
of radius $\delta$ in the complex plane provided by Fact \ref{fac:disks},
whose union contains the closed disk $\left\{ z\in\mathbb{C}:\left|z\right|\le C\right\} $.
For each $n$-tuple $\vec{\mu}=(\mu_{1},\mu_{2},....,\mu_{n})$, where
each $\mu_{j}$ is an integer in the range $1\le\mu_{j}\le M$, let
$Q(\mu_{1},\mu_{2},....,\mu_{n})=Q\left(\vec{\mu}\right)$ be the
set of all functions $g:F\to\mathbb{C}$ such that $g(a_{j})\in D_{\mu_{j}}$
for $j=1,2,...,n$. Furthermore, let $B(\mu_{1},\mu_{2},....,\mu_{n})=B\left(\vec{\mu}\right)$
be the set of all $b\in B$ for which the function $g_{b}$ defined
as above is in $Q(\mu_{1},\mu_{2},....,\mu_{n})$.

The inclusion $\left\{ z\in\mathbb{C}:\left|z\right|\le C\right\} \subset\bigcup_{q=1}^{M}D_{q}$
implies that $G\subset\bigcup_{\vec{\mu}}Q\left(\vec{\mu}\right)$
which in turn gives us that $B\subset\bigcup_{\vec{\mu}}B\left(\vec{\mu}\right)$,
i.e., that $B$ is contained in the union of the $M^{n}$ sets $B\left(\vec{\mu}\right)$.

Analogously, if $\mathbb{K}=\mathbb{R}$ then we choose a different
value for $M$, namely $M=\left\lceil \frac{C}{\delta}\right\rceil $,
and we let $I_{1}$, $I_{2}$,....., $I_{M}$ be the $M$ closed intervals
of length $2\delta$ provided by Fact \ref{fac:disks}, whose union
contains the closed interval $[-C,C]$. For each $n$-tuple $\vec{\mu}=(\mu_{1},\mu_{2},....,\mu_{n})$,
where each $\mu_{j}$ is an integer in the range $1\le\mu_{j}\le M$,
let $Q(\mu_{1},\mu_{2},....,\mu_{n})=Q\left(\vec{\mu}\right)$ be
the set of all functions $g:F\to\mathbb{\mathbb{R}}$ such that $g(a_{j})\in I_{\mu_{j}}$
for $j=1,2,...,n$. Furthermore, let $B(\mu_{1},\mu_{2},....,\mu_{n})=B\left(\vec{\mu}\right)$
be the set of all $b\in B$ for which the function $g_{b}$ defined
as above is in $Q(\mu_{1},\mu_{2},....,\mu_{n})$. As before, this
time using the inclusion $[-C,C]\subset\bigcup_{q=1}^{M}I_{q}$, we
show that $B$ is contained in the union of the $M^{n}$ sets $B\left(\vec{\mu}\right)$,
this time of course for this different choice of $M$.

Now we shall estimate the diameter of each set $B\left(\vec{\mu}\right)$.
We will use the same argument for both of the cases $\mathbb{K}=\mathbb{C}$
and $\mathbb{K}=\mathbb{R}$. Fix some $\vec{\mu}=(\mu_{1},\mu_{2},....,\mu_{n})$
and let $b$ and $b'$ be two arbitrary elements of $B\left(\vec{\mu}\right)$.
Fix some arbitrary positive number $\rho$ and let $a$ be an element
of $A$ (possibly depending on $\rho$) for which \[
d_{B}(b,b')\le\left|h(a,b)-h(a,b')\right|+\rho.\]
Next, using (\ref{eq:nfl}), we pick some $a_{j}\in F$ (possibly
depending on $a$ and therefore on $\rho$) for which $d_{A}(a,a_{j})\le\epsilon$.
Then we have \begin{eqnarray*}
\left|h(a,b)-h(a,b')\right| & \le & \left|h(a,b)-h(a_{j},b)\right|+\left|h(a_{j},b)-h(a_{j},b')\right|+\left|h(a_{j},b')-h(a,b')\right|\\
 & \le & \epsilon+\left|g_{b}(a_{j})-g_{b'}(a_{j})\right|+\epsilon\,.\end{eqnarray*}
By definition, when $\mathbb{K}=\mathbb{C}$ the two numbers $g_{b}(a_{j})$
and $g_{b'}(a_{j})$ are either both in the same disk $D_{\mu_{j}}$
of radius $\delta$, and when $\mathbb{K}=\mathbb{R}$ they are both
in the same interval $I_{\mu_{j}}$ of length $2\delta$. So, for
both $\mathbb{K}=\mathbb{C}$ and $\mathbb{K}=\mathbb{R}$, we have
\[
\left|g_{b}(a_{j})-g_{b'}(a_{j})\right|\le2\delta\]
and the previous three displayed estimates combine to give us that
$d_{B}(b,b')\le2\epsilon+2\delta+\rho$. Since $b$ and $b'$ are
arbitrary and $\rho$ can be chosen arbitrarily small, we obtain that
the diameter of $B\left(\vec{\mu}\right)$ does not exceed $2\epsilon+2\delta$
and therefore $N_{B}^{\Delta}(\epsilon+\delta)\le M^{n}$, which is
exactly (\ref{eq:cgkt}) when $\mathbb{K}=\mathbb{C}$ or exactly
(\ref{eq:gkt}) when $\mathbb{K}=\mathbb{R}$. $\qed$

\section{\label{sec:example}Two more examples}

The two examples which we present here apply to the rather more general
context of Theorem \ref{thm:elq}, rather than to Corollary \ref{cor:qschauder}.

\subsection{The estimate (\ref{eq:gkt}) in Theorem \ref{thm:elq} is best possible.\protect \\
}

Our example in this subsection is a simpler version of the example
discussed in Section \ref{sec:intro}. It shows (as also follows indirectly
from the example of Section \ref{sec:intro}) that (\ref{eq:gkt})
cannot be improved:

Let $m$ and $n$ be positive integers. Let $A$ be the set of $n$
canonical vectors $\vec{e}_{j}$, $j=1,2,...,n$ in $\mathbb{R}^{n}$,
i.e., $\vec{e}_{1}=(1,0,0,....,0)$, $\vec{e}_{2}=(0,1,0,.....,0)$,
....., $\vec{e}_{n}=(0,0,0,....,0,0,1)$. Divide the cube $[-1,1]^{n}$
in $\mathbb{R}^{n}$ into $m^{n}$ non overlapping closed cubes, each
of side length $2/m$. Let $B$ be the set of $m^{n}$ points $\vec{b}=(\beta_{1},\beta_{2},....,\beta_{n})$
which are the centres of these cubes. Let $h(\vec{a},\vec{b})$ be
the usual inner product of vectors in $\mathbb{R}^{n}$. Thus $d_{B}(\vec{b},\vec{b'})=\left\Vert \vec{b}-\vec{b'}\right\Vert _{\ell_{n}^{\infty}}$
for each $\vec{b}$ and $\vec{b}'$ in $B$, and $d_{A}(\vec{e}_{j},\vec{e}_{k})=2-2/m$
whenever $j\ne k$. This means that $N_{A}(\rho)=n$ for all $\rho\in(0,1-1/m)$.
We also have $\sup_{\vec{a}\in A,\vec{b}\in B}\left|h(\vec{a},\vec{b})\right|=1-1/m$.
Since any set in $\mathbb{R}^{n}$ whose $\ell_{n}^{\infty}$ diameter
is less than $2/m$ can contain at most one of the points of $B$,
we see that $N_{B}^{\Delta}(\rho)=m^{n}$ for each $\rho\in(0,1/m)$.

Suppose that $m\ge2$, and choose positive numbers $\epsilon=\frac{1}{3m^{2}}$
and $\delta=\frac{1}{m}-\frac{1}{2m^{2}}$. Note that $\epsilon+\delta<1/m$
and $\epsilon<1-1/m$. So, in this case, $N_{B}^{\Delta}(\epsilon+\delta)=m^{n}$
and $N_{A}(\epsilon)=n$. Furthermore the number \[
\frac{C}{\delta}=\frac{1-\frac{1}{m}}{\frac{1}{m}-\frac{1}{2m^{2}}}=\frac{m-1}{1-\frac{1}{2m}}\]
 lies in the interval $(m-1,m]$. Consequently, $\left(\left\lceil \frac{C}{\delta}\right\rceil \right)^{N_{A}(\epsilon)}=m^{n}=N_{B}^{\Delta}(\epsilon+\delta)$
and in this case equality holds in (\ref{eq:gkt}).\bigskip{}

\subsection{The statement (ii) in Theorem \ref{thm:elq} is very close to best
possible. \protect \\
}

Our example in this subsection addresses the statement (ii) in Theorem
\ref{thm:elq}. It is natural to ask whether (ii) might remain true
if we replace $\epsilon$ in the requirement $\rho>\epsilon$ by some
smaller number. In order to show that (ii) cannot be sharpened in
such a way, we will give an example of sets $A$ and $B$ and a function
$h:A\times B\to\mathbb{R}$ satisfying the hypotheses of Theorem \ref{thm:elq},
such that, when $\epsilon=1/2$, we have $N_{A}(\epsilon)<\infty$,
but $N_{B}^{\Delta}(\rho)=\infty$ for all $\rho\in(0,\epsilon)$.
(Since we have not yet determined whether or not $N_{B}^{\Delta}(1/2)$
is finite, this does not settle the question of whether or not $N_{A}(\epsilon)<\infty$
always implies that $N_{B}^{\Delta}(\epsilon)<\infty$.)

In view of (\ref{eq:tkcn}), statement (ii) of Theorem \ref{thm:elq}
implies the following similar statement, expressed solely in terms
of intrinsic covering numbers:\begin{equation}
\mbox{If }N_{A}\left(\epsilon\right)<\infty\,\mbox{then }N_{B}(\rho)<\infty\,\mbox{for all }\rho>2\epsilon\,.\label{eq:icnc}\end{equation}

The same example which we are about to present will also show that
the number $2\epsilon$ in (\ref{eq:icnc}) cannot be replaced by
an smaller number. I.e., we can have $N_{A}(\epsilon)<\infty$ but
$N_{B}(\rho)=\infty$ for all $\rho\in(0,2\epsilon)$. Here again
we will be choosing $\epsilon=1/2$. (Since in our example $N_{B}(1)=1$,
this does not settle the question of whether or not we can always
replace $\rho>2\epsilon$ by $\rho\ge2\epsilon$ in (\ref{eq:icnc}).)

Here are the details of the example:

We take $A$ to be the set of all real valued sequences $a=\left\{ \alpha_{n}\right\} _{n\in\mathbb{N}}$
which take values in $[0,1]$. We take $B$ to be the set of all real
valued sequences $b=\left\{ \beta_{n}\right\} _{n\in\mathbb{N}}$
which are finitely supported and satisfy $\sum_{n=1}^{\infty}\left|\beta_{n}\right|\le1$.
For each $a=\left\{ \alpha_{n}\right\} _{n\in\mathbb{N}}$ in $A$
and each $b=\left\{ \beta_{n}\right\} _{n\in\mathbb{N}}$ in $B$,
let $h(a,b)=\sum_{n=1}^{\infty}\alpha_{n}\beta_{n}$. Then of course
$d_{A}(a-a')=\left\Vert a-a'\right\Vert _{\ell^{\infty}}$ for all
$a,a'\in A$, and, for each $b=\left\{ \beta_{n}\right\} _{n\in\mathbb{N}}$
and $b'=\left\{ \beta'_{n}\right\} _{n\in\mathbb{N}}$, we have $d_{B}(b,b')=\max\left\{ \sum_{n=1}^{\infty}\left(\beta_{n}-\beta_{n}'\right)_{+},\sum_{n=1}^{\infty}\left(\beta_{n}'-\beta_{n}\right)_{+}\right\} $.

Let $x$ be the sequence $x=\left\{ \xi_{n}\right\} _{n\in\mathbb{N}}$
where $\xi_{n}=1/2$ for all $n$. Then $d_{A}(a,x)\le1/2$ for all
$a\in A$. This means that $N_{A}(1/2)=1<\infty$. Since $d_{B}(b,b')\le2$
for all $b,b'\in B$, we see that $N_{B}^{\Delta}(1)=1$. Since $d_{B}(b,0)\le1$
for all $b\in B$, we also see that $N_{B}(1)=1$.

Now let us show that $N_{B}(\rho)=\infty$ whenever $0<\rho<1$: Suppose
on the contrary that $N_{B}(\rho)=m<\infty$ for some $\rho\in(0,1)$.
Then there exist $m$ sequences in $B$, which we will denote by $b_{j}=\left\{ \beta_{j,n}\right\} _{n\in\mathbb{N}}$
for $j=1,2,....,m$, such that \begin{equation}
\min\left\{ d_{B}(b,b_{j}):1\le j\le m\right\} \le\rho\,\,\mbox{for all }b\in B.\label{eq:whc}\end{equation}
Since each of the sequences $b_{j}$ is finitely supported, there
exists an integer $q$ such that $\beta_{j,q}=0$ for $j=1,2,...,m$.
Now let $b=\left\{ \beta_{n}\right\} _{n\in\mathbb{N}}$ be the sequence
for which $\beta_{q}=1$ and $\beta_{n}=0$ for all $n\ne q$. Then
$d_{B}(b,b_{j})\ge1$ for all $j$, which contradicts (\ref{eq:whc}).
Thus indeed, for each $\rho\in(0,1)$, we have $N_{B}(\rho)=\infty$,
and consequently (cf.~(\ref{eq:tkcn})) also $N_{B}^{\Delta}(\rho/2)=\infty$.
In other words we also have $N_{B}^{\Delta}(\rho)=\infty$ for each
$\rho\in(0,1/2)$. This shows that the covering numbers $N_{A}$,
$N_{B}^{\Delta}$ and $N_{B}$ have all the properties which were
promised at the beginning of this subsection.

As indicated above, we would like to know whether or not $N_{B}^{\Delta}(1/2)<\infty$.

\section{\label{sec:appendix}Appendix: Further details about the case where
$A$ is absolutely convex}

Here is a proof of the estimate (\ref{eq:zefc}) which was mentioned
in Remark \ref{rem:absconv}:

Let $\epsilon$ be the positive number appearing in (\ref{eq:zefc}),
let $n=N_{A}(\epsilon)$, and let $F=\left\{ a_{1},a_{2},.....,a_{n}\right\} $
be the finite subset of $A$ depending on $\epsilon$ which is introduced
at the beginning of the proof of Theorem \ref{thm:elq}. If $\Gamma_{j}=\left\{ a\in A:d_{A}(a,a_{j})\le\epsilon\right\} $
for each $j$, then (\ref{eq:nfl}) can be restated as $A\subset\bigcup_{j=1}^{n}\Gamma_{j}$.

Since $h(0,b)=0$ for all $b\in B$ and since $d_{A}$ is related
to the seminorm $p$ by (\ref{eq:ffda}), we see that \begin{equation}
C:=\sup_{a\in A,b\in B}\left|h(a,b)\right|=\sup_{a\in A,b\in B}\left|h(a,b)-h(0,b)\right|=\sup_{a\in A}d_{A}(a,0)=\sup_{a\in A}p(a)\,.\label{eq:qcpt}\end{equation}

Let $a_{*}$ be an arbitrary fixed element of $A$. The set $L=\{(1-2t)a_{*}:t\in[0,1]\}$,
i.e., the line segment in $V$ from $a_{*}$ to $-a_{*}$, is contained
in $A$ and therefore also in $\bigcup_{j=1}^{n}\Gamma_{j}$. So the
interval $[0,1]$ satisfies \begin{equation}
[0,1]\subset\bigcup_{j=1}^{n}I_{j}\,,\label{eq:mdpt}\end{equation}
 where the sets $I_{j}$ are defined by \[
I_{j}=\{t\in[0,1]:(1-2t)a_{*}\in\Gamma_{j}\}=\left\{ t\in[0,1]:p\left(a_{*}-2ta_{*}-a_{j}\right)\le\epsilon\right\} \,.\]
For each $j$ the non negative function $t\mapsto p\left(a_{*}-2ta_{*}-a_{j}\right)$
is continuous, and in fact it is also convex. So each $I_{j}$ is
a closed (possibly empty) subinterval of $[0,1]$. It follows from
(\ref{eq:mdpt}) that at least one of these intervals, say $I_{m}$,
must have length no less than $1/n$. Let $s$ be a point in $I_{m}$
such that $s+1/n$ is also in $I_{m}$. Then \begin{eqnarray*}
\frac{2}{n}p(a_{*}) & = & p\left(\frac{2a_{*}}{n}\right)=p\left(\left(a_{*}-2sa_{*}-a_{m}\right)-\left(a_{*}-2\left(s+\frac{1}{n}\right)a_{*}-a_{m}\right)\right)\\
 & \le & p\left(a_{*}-2sa_{*}-a_{m}\right)+p\left(a_{*}-2\left(s+\frac{1}{n}\right)a_{*}-a_{m}\right)\le\epsilon+\epsilon\,.\end{eqnarray*}
This shows that $p(a_{*})\le\epsilon n$. So now we simply take the
supremum in this last inequality as $a_{*}$ ranges over $A$ and
use (\ref{eq:qcpt}) to obtain (\ref{eq:zefc}).

\medskip{}

\textbf{\textit{Acknowledgement.}} We are most grateful for some helpful
comments from Emanuel and Vitali Milman about the previous version
of this paper, and also from Mario Milman during an earlier stage
of this research.

\end{document}